# THE FOOTBALL PLAYER AND THE INFINITE SERIES

HAROLD P. BOAS

## 1. THE FOOTBALL PLAYER

The air buzzed with anticipation as the football team crowded excitedly into the lecture hall. The country's top halfback was about to defend his PhD thesis in mathematics! It soon became apparent that the proceedings were a mere formality, as the candidate's dissertation on summability methods for divergent Dirichlet series was a masterful piece of work.

This scenario is no fantasy from a 1990's television sitcom: it is a true story. The place was Copenhagen, the year was 1910, and the sport was "football" as the word is understood internationally ("soccer" in American lingo). The star halfback played in the 1908 Olympics on Denmark's silver-medal football team, a team that is still in the record books [21, p. 172] for the most goals scored in a single game. (Denmark defeated France by the lopsided score of 17 to 1.) The dissertation title was *Contributions to the Theory of Dirichlet Series* (well, actually *Bidrag til de Dirichlet'ske Rækkers Theori*), and the candidate's name was Harald Bohr.

(Devotees of American football remember Frank Ryan, who wrote his PhD dissertation [23, 24] on geometric function theory while quarterback for the Cleveland Browns, champions of the National Football League at the time. But that's another story [18, 22].)

Among mathematicians, Harald Bohr is best remembered today for his theory of almost periodic functions [10]; students of complex analysis also know him for the Bohr-Mollerup theorem (see, for example, [3, Theorem 2.1], [12, §§274–275]) that characterizes the $\Gamma$ function on the positive real axis as the unique positive, logarithmically convex function $f$ such that $f(x+1) = xf(x)$ for all $x$ and $f(1) = 1$. In his native land, Bohr's early fame as a sports hero and his subsequent prominence as a distinguished academician were eclipsed by his status as the kid brother of Niels Bohr. Brother Niels, a prime architect of

Research partially supported by NSF grant number DMS 9500916.

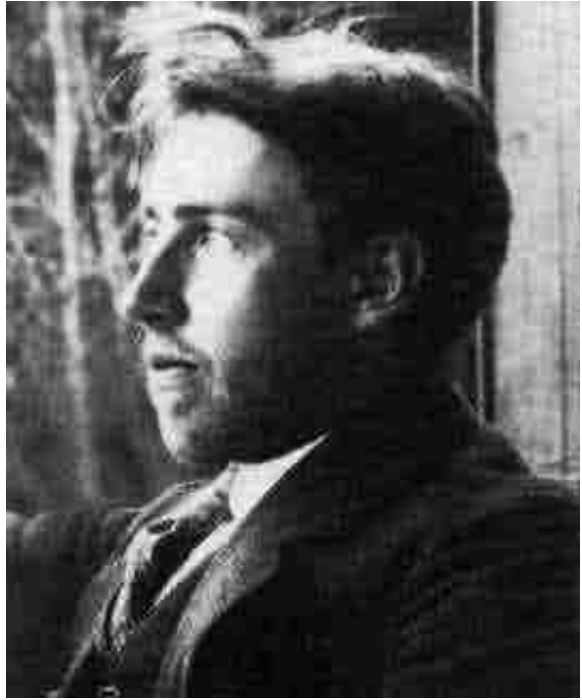

Figure 1. Harald Bohr

modern atomic theory and recipient of the Nobel prize for physics in 1922, was Denmark's most honored citizen during his lifetime.

## 2. The infinite series

Like many others before and after him, Harald Bohr wanted to decide the truth or falsity of the Riemann hypothesis, one of the most famous unsolved problems of mathematics. Bohr was unsuccessful, but much of his mathematical work was motivated by trying to understand the Riemann zeta-function $\zeta$:

$$\zeta(s) = \sum_{n=1}^{\infty} \frac{1}{n^s}, \qquad \operatorname{Re} s > 1.$$

It is easy to see that the infinite series on the right-hand side converges absolutely in the half-plane where the real part of the complex variable $s$ exceeds 1, for $|1/n^s| = 1/n^{\operatorname{Re} s}$, and $\sum_{n=1}^{\infty} 1/n^x$ converges when $x > 1$. On the other hand, there is no larger open half-plane where the series converges (even conditionally), because when $s = 1$ the series reduces to the divergent harmonic series.

It is a natural idea to try to understand the Riemann $\zeta$-function by studying the more general *Dirichlet series* of the form $\sum_{n=1}^{\infty} a_n/n^s$, the coefficients $a_n$ being complex constants. (These are *ordinary* Dirichlet series; for a wider class, see, for example, [2, 16].) A simple example of a Dirichlet series is $\sum_{n=1}^{\infty}(-1)^{n+1}/n^s$, which is the $\zeta$-function series with alternating signs. Evidently this series converges absolutely in exactly the same half-plane as the $\zeta$-function series does: $\operatorname{Re} s > 1$.

However, this new series converges conditionally (but not absolutely) in the larger half-plane where $\operatorname{Re} s > 0$. The convergence follows from the Abel-Dirichlet-Dedekind generalization of the alternating series test (see, for example, [15, §143], [20, §5.5]), which implies that if $\{b_n\}$ is a sequence tending to 0 and of bounded variation, then $\sum_n (-1)^n b_n$ converges. (The sequence $\{1/n^s\}$ has bounded variation because $|1/n^s - 1/(n+1)^s| = O(1/n^{1+\operatorname{Re} s})$, and $\sum_n 1/n^{1+\operatorname{Re} s}$ converges when $\operatorname{Re} s > 0$.)

This phenomenon of conditional convergence is contrary to our experience with ordinary power series $\sum_{n=1}^{\infty} c_n z^n$, for a power series converges absolutely at all points of its open disk of convergence. A Dirichlet series can converge nonabsolutely (that is, conditionally) in a vertical strip, and the above example shows that the width of such a strip can be as large as 1. The width of the strip of conditional, nonabsolute convergence cannot, however, exceed 1. Indeed, if $\sum_{n=1}^{\infty} a_n/n^s$ converges for a certain $s$, then the individual terms tend to 0, and in particular are bounded in absolute value by some constant $M$; now if $z$ is a complex number such that $\operatorname{Re} z > 1 + \operatorname{Re} s$, then $\sum_{n=1}^{\infty} |a_n/n^z| \leq M \sum_{n=1}^{\infty} 1/n^{\operatorname{Re}(z-s)} < \infty$.

It happens that the series $\sum_{n=1}^{\infty}(-1)^{n+1}/n^s$ is closely related to the $\zeta$-function. When $\operatorname{Re} s > 1$, we can rearrange the terms of the absolutely convergent series however we like, so by separating the sum over odd integers from the sum over even integers, we find that

$$\sum_{n=1}^{\infty} \frac{(-1)^{n+1}}{n^s} = \sum_{n=1}^{\infty} \frac{1}{n^s} - 2\sum_{k=1}^{\infty} \frac{1}{(2k)^s} = \zeta(s)(1 - 2^{1-s}).$$

Thus, the function $(1 - 2^{1-s})^{-1} \sum_{n=1}^{\infty}(-1)^{n+1}/n^s$ serves to extend the definition of the $\zeta$-function from the half-plane where $\operatorname{Re} s > 1$ to the half-plane where $\operatorname{Re} s > 0$. Bohr observed in [6] that one way to extend the $\zeta$-function to the whole plane is to take iterated Cesàro averages of the series $\sum_{n=1}^{\infty}(-1)^{n+1}/n^s$. The famous Riemann hypothesis is equivalent to the statement that the zeroes of the function $\sum_{n=1}^{\infty}(-1)^{n+1}/n^s$ in the right half-plane all lie on the vertical line where $\operatorname{Re} s = 1/2$.

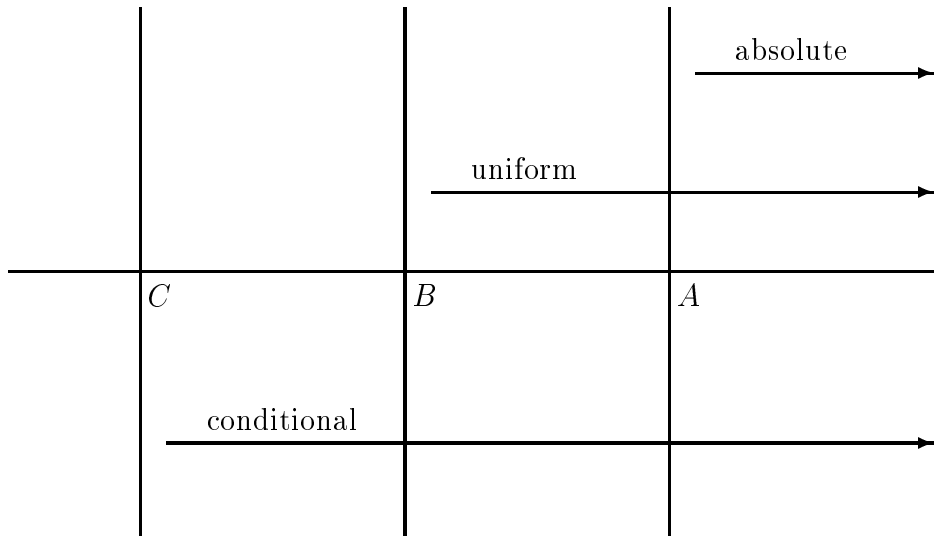

FIGURE 2. Convergence regions for Dirichlet series

## 3. THE QUESTION

What about uniform convergence of Dirichlet series? An ordinary power series converges uniformly on each closed disk inside its open disk of convergence, but this gives no hint about what might be true for Dirichlet series (as we have already seen in the case of conditional convergence).

Since $|1/n^s|$ does not depend on the imaginary part of $s$, it is clear that if a Dirichlet series $\sum_{n=1}^{\infty} a_n/n^s$ converges absolutely in a half-plane where $\operatorname{Re} s > A$, then it converges uniformly in each closed half-plane $\{s : \operatorname{Re} s \geq A + \epsilon\}$, where $\epsilon$ can be any positive number. Having just seen that there may be an abscissa $C$ to the left of $A$ such that the series converges conditionally when $\operatorname{Re} s > C$, we might anticipate that there is an intermediate abscissa $B$, as indicated in Figure 2, such that the Dirichlet series converges uniformly in each closed half-plane $\{s : \operatorname{Re} s \geq B + \epsilon\}$, where $\epsilon > 0$. Harald Bohr introduced this notion of a line of uniform convergence in [7].

In [8, p. 446], Bohr asked: what is the maximal possible width $A - B$ of the vertical strip of uniform, but not absolute, convergence of a Dirichlet series? We saw above that $A - C \leq 1$, so certainly $A - B$ cannot exceed 1. It turns out that $A - B$ cannot exceed $1/2$, and this value is sharp.

Although Bohr knew that $A - B \leq 1/2$, he could not produce a single example for which $A - B > 0$. In a companion paper [26] in the same volume, Otto Toeplitz gave examples showing that the upper cut-off for $A - B$ is no smaller than $1/4$. It was nearly two decades later that H. F. Bohnenblust and Einar Hille finally proved $1/2$ to be the right value in an article [5] that Henry Helson, commenting in Hille's collected papers [19, p. 664], termed "a remarkable piece of work." The result was rediscovered by Seán Dineen and Richard M. Timoney more than half a century later [14], with a new proof based on the relationship [13] between nuclearity and the existence of absolute bases in certain locally convex spaces.

My aim here is to make the theorem accessible to a wide audience by presenting a relatively elementary proof using only methods that have existed in textbook form since I was in high school. (I hesitate to call the methods "classical," however. The technique of random polynomials discussed below in paragraph 5.1.3 was not available to Bohr.) My attention was directed to this theorem by Henry Helson when I lectured at Berkeley about some joint work with Dmitry Khavinson [4] concerning another problem of Bohr.

## 4. The upper bound

Suppose that a Dirichlet series $\sum_{n=1}^{\infty} a_n/n^s$ converges uniformly on a vertical line where $\operatorname{Re} s = b$. I claim that if $\epsilon$ is an arbitrary positive number, then the series converges absolutely when $\operatorname{Re} s \geq b + \epsilon + 1/2$; that is, $\sum_{n=1}^{\infty} |a_n|/n^{b+\epsilon+1/2} < \infty$. In other words, the width $A - B$ is no larger than $1/2$.

Observe that by the Cauchy-Schwarz inequality, $\sum_{n=1}^{\infty} |a_n|/n^{b+\epsilon+1/2}$ is at most $(\sum_{n=1}^{\infty} |a_n|^2/n^{2b})^{1/2} (\sum_{n=1}^{\infty} 1/n^{1+2\epsilon})^{1/2}$. Since $\sum_{n=1}^{\infty} 1/n^{1+2\epsilon}$ converges, the claim will follow if I show that $\sum_{n=1}^{\infty} |a_n|^2/n^{2b}$ converges.

Since each finite partial sum $\sum_{n=1}^{N} a_n/n^s$ is bounded on the line where $\operatorname{Re} s = b$, and since (by hypothesis) the partial sums converge uniformly on this line, the partial sums must be uniformly bounded on the line, say by a constant $M$. Then for every positive integer $N$ and every real number $t$, we have the inequality

$$M^2 \geq \left| \sum_{n=1}^{N} \frac{a_n}{n^{b+it}} \right|^2 = \sum_{n=1}^{N} \frac{|a_n|^2}{n^{2b}} + 2 \operatorname{Re} \sum_{1 \leq n < m \leq N} \frac{a_n \bar{a}_m}{(nm)^b (n/m)^{it}}.$$

Taking the average value with respect to $t$ by integrating from $-T$ to $T$ and dividing by $2T$, we find that

$$M^2 \geq \sum_{n=1}^{N} \frac{|a_n|^2}{n^{2b}} + 2 \operatorname{Re} \sum_{1 \leq n < m \leq N} \frac{a_n \bar{a}_m}{(nm)^b} \frac{\sin(T \log(m/n))}{T \log(m/n)}.$$

Taking the limit as $T \to \infty$ shows that $M^2 \geq \sum_{n=1}^{N} |a_n|^2/n^{2b}$. Since $N$ is arbitrary, this means that $\sum_{n=1}^{\infty} |a_n|^2/n^{2b}$ does converge.

This confirms that the maximal width $A - B$ of the strip of uniform but not absolute convergence of a Dirichlet series is at most $1/2$. Next I want to show that the cut-off value for this width is no smaller than $1/2$.

## 5. The lower bound

I will construct a Dirichlet series $\sum_{n=1}^{\infty} a_n/n^s$ that converges uniformly in every half-plane $\{s : \operatorname{Re} s \geq \delta + \frac{1}{2}\}$, where $\delta > 0$, but that does not converge absolutely when $\operatorname{Re} s < 1$. This example demonstrates that no number smaller than $1/2$ will serve as a cut-off for the maximal width of the strip of uniform nonabsolute convergence of Dirichlet series.

5.1. **Tools.** The construction uses off-the-shelf technology: elementary counting, the prime number theorem, and the theory of random Fourier series. There is enough slack in the method that I do not need particularly sharp implementations of these tools. The theory of analytic functions of an infinite number of variables, central to Harald Bohr's approach, is hiding in the background, but I shall not need to make explicit reference to it.

Nonetheless, the philosophy of the construction is very much that of Bohr. Namely, I choose to view an object such as $1/45^s$ not as the reciprocal of a power of an integer, but as the value of the monomial $z_1^2 z_2$ when $z_1 = 1/3^s$ and $z_2 = 1/5^s$. Thus, the problem becomes separated from number theory and turns into a problem about polynomials.

5.1.1. *The prime number theorem.* The most familiar version of the prime number theorem says that the number of primes less than $x$ is asymptotic to $x/\log x$ when $x \to \infty$. An equivalent statement is that if the prime numbers are arranged in increasing order ($p_1 = 2$, $p_2 = 3$, $p_3 = 5$, and so on), then the size of the $n$th prime $p_n$ is asymptotic to $n \log n$. I need only the weaker statement that there is a constant $c_1$ larger than 1 such that $1/c_1 < p_n/(n \log n) < c_1$ when $n > 1$, which is rather easier to prove than the full-blown prime number theorem (see, for example, [1, §4.5]).

5.1.2. *Counting monomials.* I need simple bounds on the number of monomials of degree $m$ in $n$ variables: objects of the form $z_1^{\alpha_1} z_2^{\alpha_2} \ldots z_n^{\alpha_n}$, where the $\alpha_j$ are nonnegative integers whose sum is $m$. Viewing such a monomial as a product of $m$ nontrivial factors, where there are $n$ choices for each factor, gives a count of $n^m$; but this count is too big, since it takes account of the order of the terms. No particular product of terms has more than $m!$ rearrangements, and some products have fewer rearrangements, so $n^m/m!$ is an undercount. Thus the number of distinct monomials of degree $m$ in $n$ variables is between $n^m/m!$ and $n^m$. It is easy to show that the precise count is the binomial coefficient $\binom{n+m-1}{m}$, but I shall not need this exact value.

5.1.3. *Random polynomials.* Consider a homogeneous polynomial of degree $m$ in $n$ complex variables with coefficients $\pm 1$: that is, an object of the form
$$\sum_{\alpha_1+\alpha_2+\cdots+\alpha_n=m} \pm z_1^{\alpha_1} z_2^{\alpha_2} \ldots z_n^{\alpha_n}.$$
To avoid trivialities, I assume that $m$ and $n$ are both at least 2. What can be said about the supremum of the modulus of such a polynomial when every coordinate $z_j$ lies in the unit disk? Since each term has modulus at most 1, the maximum modulus is certainly no more than the total number of terms, which according to the previous paragraph is less than $n^m$. However, such a polynomial typically has maximum modulus much smaller than this crude bound. According to the theory of random trigonometric polynomials (see, for example, [17, Theorem 4 of Chapter 6]), there is a constant $c_2$ such that if the $\pm$ signs are assigned at random, then with high probability the maximum modulus of the resulting polynomial is less than $c_2 n^{(m+1)/2} \sqrt{\log m}$. Although there consequently are many polynomials satisfying this bound, all I need is the existence of one for each $m$ and $n$.

5.2. **The construction.** I will construct a Dirichlet series $\sum_{n=1}^{\infty} a_n/n^s$ for which every coefficient $a_n$ is either $0$, $+1$, or $-1$, and I will show that this Dirichlet series converges uniformly when $\operatorname{Re} s \geq \delta + \frac{1}{2}$, where $\delta > 0$, yet does not converge absolutely when $\operatorname{Re} s < 1$.

I construct the terms of the series in groups. To build the $k$th group (starting with $k = 2$), choose a random homogeneous polynomial of degree $k$ in $2^k$ variables with coefficients $\pm 1$ (as described in paragraph 5.1.3). List the $2^k$ consecutive prime numbers starting with the $2^k$th prime, and for each such prime $p$, substitute $1/p^s$ for the corresponding variable in the polynomial. This converts the sum of monomials $\pm z_1^{\alpha_1} z_2^{\alpha_2} \ldots z_{2^k}^{\alpha_{2^k}}$ into a sum of terms $\pm 1/n^s$, where each integer $n$ is the product of exactly $k$ primes (counting repeated factors with their multiplicities) from the block of $2^k$ primes starting at the $2^k$th prime. The uniqueness of prime factorization implies that no integer $n$ appears more than once.

For every integer $n$ not arising in the above process, I set $a_n = 0$. The first integer $n$ for which $a_n \neq 0$ is 49, for this is the smallest integer that is the product of 2 primes taken from the set of $2^2$ consecutive primes starting with $p_4 = 7$.

Now I verify that the constructed Dirichlet series has the required properties. First consider the question of absolute convergence of $\sum_{n=49}^{\infty} a_n/n^s$. The counting argument in paragraph 5.1.2 implies that the number of integers $n$ formed from products of $k$ primes in the block from the $2^k$th prime to the $2^{k+1}$th prime exceeds $2^{k^2}/k^k$. By the prime number theorem, the $2^{k+1}$th prime is bounded above by $3c_1 k 2^k$, so such integers $n$ are bounded above by $(3c_1 k)^k 2^{k^2}$. Hence $\sum_{n=49}^{\infty} |a_n/n^s|$ exceeds $\sum_{k=2}^{\infty} 2^{k^2(1-\operatorname{Re} s)}/(3c_1 k)^{k(1+\operatorname{Re} s)}$. Evidently the latter sum diverges when $\operatorname{Re} s < 1$, so our Dirichlet series fails to converge absolutely when $\operatorname{Re} s < 1$. (On the other hand, since the coefficients $a_n$ are bounded, it is evident that our Dirichlet series does converge absolutely when $\operatorname{Re} s > 1$.)

Next consider the question of uniform convergence of our Dirichlet series. I wish to estimate the modulus of the sum of the terms in the $k$th block. This piece of the Dirichlet series equals the value of our random polynomial when we substitute for each variable the reciprocal of the corresponding prime number raised to the power $s$. Since the polynomial is homogeneous of degree $k$, the supremum of its modulus when the variables have modulus at most $|1/p^s|$ is $1/p^{k \operatorname{Re} s}$ times the bound $c_2 2^{k(k+1)/2} \sqrt{\log k}$ coming from paragraph 5.1.3. Since the $2^k$th prime is bounded below by $k 2^k/2c_1$, this chunk of the Dirichlet series is bounded above by $c_2 2^{k(k+1)/2} \sqrt{\log k}/(k 2^k/2c_1)^{k \operatorname{Re} s}$. The Weierstrass $M$-test and the root test now imply that the series of blocks converges uniformly when $\operatorname{Re} s \geq 1/2$.

The proof is now finished modulo a technical (but nontrivial) point. I have showed that the constructed Dirichlet series converges uniformly for $\operatorname{Re} s \geq 1/2$ if the series is summed in appropriate blocks; however, I need to show that the Dirichlet series converges uniformly when summed in its natural order, without grouping. This follows from a general lemma, essentially due to Bohr [9, Hilfssatz 2].

**Lemma.** *Suppose that a Dirichlet series $\sum_{n=1}^{\infty} b_n/n^s$ converges absolutely when $\operatorname{Re} s > a$ and that the analytic function $f(s)$ which it represents continues analytically to the half-plane where $\operatorname{Re} s > c$. If $c < b < a$, and if $f$ is bounded on the half-plane where $\operatorname{Re} s \geq b$, then for every positive $\delta$, the Dirichlet series converges uniformly on the half-plane where $\operatorname{Re} s \geq b + \delta$.*

In our situation, the series summed in blocks converges uniformly in the closed half-plane where $\operatorname{Re} s \geq \frac{1}{2}$ to a bounded function $f$ that is analytic in the open half-plane. When $\operatorname{Re} s > 1$, this function $f$ does equal the sum of the Dirichlet series (summed in any order, since in that region the series converges absolutely). Consequently, the lemma implies that the Dirichlet series converges uniformly to $f$ in each half-plane where $\operatorname{Re} s \geq \delta + \frac{1}{2}$.

5.3. **Proof of the lemma.** The lemma follows from a technique that Bohr attributed to his contemporary W. Schnee, who wrote his dissertation in Berlin in 1908 under the influence (although not the formal tutelage) of the famous Edmund Landau. Soon after receiving his Master's degree in 1909, Bohr himself began a collaboration with Landau, who had just been appointed Minkowski's successor at the University of Göttingen. In his reminiscences [11, p. xxvi], Bohr remarked on Landau's unexcelled zeal:

> When Landau and I thought that an oral conference on our work was needed, I caught the train to Göttingen for a few days' stay. No one could be in such an excellent mood for work as Landau, and his speed and perseverance were sometimes quite breathtaking. In order to show me at once that the time had come for serious work, he had instituted the tradition of ringing the bell immediately, as soon as I had arrived at his house after the long and somewhat tiring journey and had set foot inside his study, and of requesting the entering maid to inform the kitchen that 'tonight at 2 AM a very strong cup of coffee is to be served to both of us.'

The idea of the proof is easier to describe than to implement: integrate over a vertical contour, and use Cauchy's integral formula to push the contour to the right into the region where the Dirichlet series is already known to converge uniformly. The technique is still the standard one employed to derive Perron's formula for the partial sums of Dirichlet series (see, for example, [1, §11.12], [25, §9.42]). This shows that contour integration remains useful, even though symbolic computation software packages such as Mathematica® and Maple® can now calculate all the real integrals that are given in complex analysis textbooks as the main applications of contour integration.

To begin the proof, let $K$ denote an upper bound for $f$ in the half-plane where $\operatorname{Re} s \geq b$, and fix a positive $\delta$ (which we may as well assume is less than 1). I aim to show that if $\operatorname{Re} s \geq b + \delta$, then $|f(s) - \sum_{n=1}^{M} b_n/n^s|$ is bounded by a constant times $M^{-\delta} \log M$, where the constant depends on $K$ and $\delta$, but is independent of $s$ and $M$. Consequently, the Dirichlet series will converge uniformly to $f$ in the half-plane where $\operatorname{Re} s \geq b + \delta$, as claimed.

Viewing $s$ and $M$ as fixed for the moment, with $\operatorname{Re} s \geq b + \delta$, consider integrating $f(z)(M+\frac{1}{2})^{z-s}/(z-s)$ as a function of $z$ around the rectangular contour shown in Figure 3 with vertices at $s-\delta-iM^{a-b+2}$, $s+a-b-iM^{a-b+2}$, $s+a-b+iM^{a-b+2}$, and $s-\delta+iM^{a-b+2}$. By

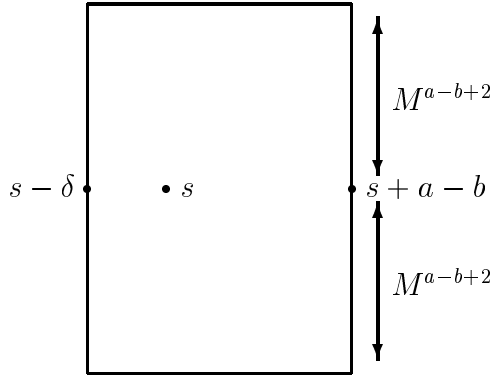

FIGURE 3. Integration contour

Cauchy's integral formula, this integral equals $2\pi i f(s)$. The integral over the left-hand edge of the rectangle has modulus bounded by $KM^{-\delta}\int_{-M^{a-b+2}}^{M^{a-b+2}}(\delta^2+y^2)^{-1/2}\,dy$, and hence by a constant (depending on $\delta$ and $K$) times $M^{-\delta}\log M$. The integrals over the top and bottom edges of the rectangle are each bounded by

$$KM^{-(a-b+2)}\int_{-\delta+\operatorname{Re} s}^{a-b+\operatorname{Re} s}(M+\tfrac{1}{2})^{x-\operatorname{Re} s}\,dx,$$

and hence by a constant times $M^{-2}$. Consequently, $2\pi i f(s)$ differs by $O(M^{-\delta}\log M)$ from the integral over the right-hand edge of the rectangle.

Since the right-hand edge of the contour is in the region where the Dirichlet series is known to converge uniformly to $f$, we may replace $f(z)$ by $\sum_{n=1}^{\infty} b_n/n^z$ in the remaining integral and interchange the order of summation and integration. We now have that

$$(1) \quad \left| 2\pi i f(s) - \sum_{n=1}^{\infty}\frac{b_n}{n^s}\int_{s+a-b-iM^{a-b+2}}^{s+a-b+iM^{a-b+2}}\left(\frac{M+\frac{1}{2}}{n}\right)^{z-s}\frac{1}{z-s}\,dz \right|$$

is $O(M^{-\delta}\log M)$. To evaluate the integrals in this sum, we must distinguish between the cases $n \geq M+1$ and $n \leq M$.

When $n \geq M+1$, build a new rectangular contour whose left-hand edge is the given vertical line segment with abscissa $a-b+\operatorname{Re} s$ and whose right-hand edge has very large abscissa. The integrand has no singularities inside this contour, so the integral over the left-hand side equals the negative of the sum of the integrals over the other three sides. Since $((M + \frac{1}{2})/n)^{z-s}$ is decaying exponentially when $\operatorname{Re} z$ becomes large, we may push the right-hand edge of the contour off to $+\infty$. The integrals over the top and bottom sides are each bounded by $M^{-(a-b+2)} \int_{a-b}^{\infty} ((M + \frac{1}{2})/n)^x \, dx$. Hence the terms for which $n \geq M + 1$ make a total contribution to the sum in (1) not exceeding twice

$$\sum_{n \geq M+1} \frac{|b_n|}{n^{\operatorname{Re} s}} M^{-(a-b+2)} \left(\frac{M + \frac{1}{2}}{n}\right)^{a-b} \left|\log \frac{M + \frac{1}{2}}{n}\right|^{-1}.$$

Observe that $|\log(M + \frac{1}{2})/n|$ is smallest when $(M + \frac{1}{2})/n$ is closest to 1, which happens when $n = M + 1$. In this case, the absolute value of the logarithm is

$$-\log \frac{2M + 1}{2M + 2} = -\log\left(1 - \frac{1}{2M + 2}\right) > \frac{1}{2M + 2}.$$

Since the series $\sum_{n=1}^{\infty} |b_n|/n^{a-b+\operatorname{Re} s}$ is uniformly bounded above by $\sum_{n=1}^{\infty} |b_n|/n^{a+\delta}$, which converges by hypothesis, it follows that the terms in (1) with $n \geq M + 1$ have a sum bounded by a constant times $1/M$.

For the terms with $n \leq M$, we may similarly build a rectangular contour whose right-end edge is the vertical line with abscissa $a - b + \operatorname{Re} s$ and which extends to the left toward $-\infty$. The integral over this contour picks up a contribution $2\pi i$ from the simple pole at $z = s$ with residue 1, while the integrals over the top and bottom edges admit estimates analogous to the previous case. Consequently,

$$\left| f(s) - \sum_{n=1}^{M} \frac{b_n}{n^s} \right| = O(M^{-\delta} \log M)$$

uniformly with respect to $s$ when $\operatorname{Re} s \geq b + \delta$. This completes the proof of the lemma.

## 6. Envoi

We have seen an example of a Dirichlet series $f(s)$ whose strip of uniform, but not absolute, convergence attains the maximal possible width of $1/2$. On the other hand, for the Riemann zeta function $\zeta$, the width of this strip is 0. Bohnenblust and Hille went to some trouble in [5, pp. 618–620] to demonstrate that if $\lambda$ is any real number between 0 and $1/2$, then there is a Dirichlet series whose strip of uniform, non-absolute convergence has width precisely $\lambda$. Harald Bohr [5, p. 622 footnote] cut through this problem with a knife: the Dirichlet series for $f(s) + \zeta(s + \lambda)$ does the job.

DEPARTMENT OF MATHEMATICS, TEXAS A&M UNIVERSITY, COLLEGE STATION, TX 77843

*E-mail address*: boas@math.tamu.edu